# Navigating Epidemic Mathematics: Exploring Tools for Mathematical Modelling in Biology

*Pabel Shahrear, Md. Shahedul Islam, Md. Abu Bakkar, Anika Bushra, Ismail Hossain*

**Department of Mathematics, Shahjalal University of Science and Technology, Sylhet Bangladesh**

## Abstract

The ever-changing world of disease study heavily relies on mathematical models. They're key in finding and controlling infectious diseases. Our project aims to explore these mathematical tools used for studying disease spread in biology. The SEIR model holds our focus. It's a super important tool known for being flexible and useful. We kick off the study with a deep look at the modified SEIR model's design and analysis. We dive right into vital parts like the equations that make the modified SEIR model work and setting parameter identities and then checking the positivity and the limits of its solutions. The study begins with a detailed examination of the design and analysis of a modified SEIR model, demonstrating its angularity. We delve into the model's heart, dealing with critical issues such as the equations that drive the modified SEIR model, establishing parameter identities, and ensuring the positivity and boundedness of its solutions. The determination of the Basic Reproduction Number marks a significant milestone by revealing data on the disease's transmission potential. The voyage continues with an investigation into the model's local stability that includes both DFE (DFE) and EE (EE). The Global stability, a paramount consideration in understanding the long-term behaviors of the systems is scrutinized employing the Lyapunov stability theorem. The project further extends its gaze to the bifurcation analysis which is classifying and elucidating the fundamental concepts therein. One-dimensional and two-dimensional bifurcations, forward and backward bifurcation analyses are intricately examined, providing a comprehensive understanding of the system's dynamical behavior along with the basic concepts. In summarizing, this project not only offers a thorough description and analysis of the SEIR model but also lays the groundwork for advancing mathematical modeling in epidemiology. By bridging theoretical insights with practical implications this study strives to empower researcher's and policymakers with a deep understanding of infectious diseases dynamics and thereby contributing to the development of more effective and targeted public health strategies.

## Introduction
A brief introduction of epidemiology, mathematical modelling and its purpose and finally the orientations of the article is depicted.

### 1.1: Epidemy and Epidemiology

Embarking on our intellectual journey, we set out to unravel the profound significance of two of the basic concepts: epidemy and epidemiology [1]. At the heart of our exploration lies the quest to provide an easy understanding of these terms by creating a solid foundation for the comprehensive examination of infectious disease dynamics that follows.

### 1.1.1: Understanding Epidemy

Epidemy is defined as the incidence of more occurrences of a certain disease within a population, geographic area or community than is typically anticipated [1]. It is an important indication of the ups and downs of diseases, acting as a warning system for possible health problems. This sections seek to demystify the phrase and clarify the implications and offer light on how epidemic detection is the foundation of effective public health responses.

### 1.1.2: Understanding Epidemiology

Epidemiology, the field of study that supports our investigation, is concerned with the distribution, determinants, and consequences of health-related occurrences in communities [1]. Epidemiology serves as a compass for navigating the complicated landscapes of public health allowing us to examine illness trends and causes. This section tries to make epidemiology more approachable by explaining its role in understanding the complexities of infectious diseases. Influential contributions in epidemiology shape our thinking as we travel this conceptual field. This foundational work guides our exploration providing insights from the pioneering of epidemiology research. Integrating these views strengthens the narrative and provides a thorough understanding of the concepts guiding our intellectual journey.

### 1.2: The Need for Mathematical Modeling in Epidemiology

Understanding the dynamics of infectious diseases demands a sophisticated methodology. This section dives into the reasons why mathematical modeling is essential in epidemiology. We negotiate the complicated rules of infectious disease dynamics, revealing the critical importance of mathematical models. These models, which include simulation, analysis, and prediction, are excellent tools for understanding the complicated patterns of disease spread within populations.

**The Purpose and Impact of Mathematical Modeling:**
As we investigate the mathematical arena of epidemiology the aim and impact of mathematical modeling become clear [2]. This part delves into the various contributions of these approaches to public health. From influencing policy and guiding resource allocation to developing efficient ways for minimizing the impact of infectious diseases, mathematical modeling emerges as a key weapon in the fight against epidemics.

### 1.3: Significance of Study

As we begin on this mathematical travel in the field of epidemiology, the significance of our findings goes beyond theoretical exploration [3]. This effort has major consequences and contributions to the larger fields of infectious disease dynamics, public health, and mathematical modeling. Here are some significant points that highlight the relevance of our endeavor:

**1. Advancing Theoretical Understanding:** Our strategy contributes to the theoretical basis of epidemiological modeling, specifically the SEIR model [4]. We hope to improve our theoretical understanding of disease dynamics by investigating local and global stability and the bifurcation analyses and other mathematical complexities. It can feed future modeling attempts and improve our understanding of infectious disease propagation.

**2. Informing Public Health Strategies:** The findings drawn from our comprehensive evaluation of the SEIR model have practical implications for public health efforts [5].

Understanding the stability and bifurcation patterns of disease dynamics can help forecast and manage outbreaks as well as develop focused interventions and allocate resources more effectively. This project is a helpful resource for public health practitioners looking for evidence-based ways to reduce the impact of infectious diseases.

**3. Guiding Policy Decisions:** As public health policies adapt to address new concerns; our research provides a mathematical structure to inform policy decisions [3]. The investigation of global stability and bifurcation analysis provides a quantitative prism through which policymakers can evaluate the efficacy of different actions. It contributes to the development of flexible and proactive strategies for preventing the spread of infectious diseases.

**4. Empowering Future Research Endeavors:** Our study lays the groundwork for future research endeavors by initiating the exploration of the mathematical complexities of the SEIR model [5]. Researchers and academics in mathematics, epidemiology, and public health can build on our discoveries to broaden our knowledge and better comprehend the intricate relationship between mathematics and infectious illnesses.

The research is significant not only because of its extensive mathematical exploration but also because of its potential to have a good impact in the real world. It specifically seeks to transform the landscape of infectious disease management and prevention.

## 1.4: Objectives

In the field of epidemiology, our endeavor is a mathematical journey to understand infectious disease dynamics [4]. There are three overarching objectives:
**1. Examine Local Stability:** Check the SEIR model's stability, with a focus on DFE (DFE) [4]. Our goal is to use local dynamics analysis to identify factors that influence disease trajectories in specific populations.
**2. Explore Global Stability:** Expand the SEIR model's vision to include global stability by analyzing both DFE and EE (EE) [4]. Using modern mathematical tools identify overarching stability principles that govern long-term infectious disease behavior beyond localized dynamics.
**3. Bifurcation Analyses and Mathematical Intricacies:** Decipher mathematical complexities by delving into bifurcation analyses and other phenomena [4]. Unravel bifurcation patterns, understand forward and backward analyses, and investigate bifurcation categorization. Aspiring to contribute to the study of various infectious illness trajectories.
These aims serve as the foundation of our mathematical endeavor, providing an organized method to decipher infectious disease dynamics using the SEIR model. Our mission is clear as we expand theoretical understanding, with insights aimed at enriching academic discourse and influencing real-world policies in the battle against infectious diseases.

## 1.5 Outline

Within the realm of mathematical modeling and biology, this project orchestrates a harmonious blend, aiming to deepen our understanding of infectious diseases [2]. Formulation and Analysis delves into Formulation and Analysis. It starts with a explores Model Formulation and treating it like composing a mathematical melody for disease transmission nuances. We examine Equations for the Modified SEIR Model in which symbols and parameters depict the disease dynamics complex choreography. Parameter Identity ensures correctness in modeling real-world scenarios. Positivity and boundedness of solutions ensure realistic outcomes which mirrors constraints in a dance. The chapter crescendos with Calculating the Basic Reproduction Number, quantifying a diseases potential like determining a play notes resonance.

Bifurcation unfolds as a captivating harmony of Bifurcation Analysis and its Classification. Beginning with an overview, the investigation transcends into Bifurcation in one and two dimensions. This dynamic chapter unfolds with the exploration of Forward Bifurcation Analysis, contributing to our quest for mathematical precision and biological insight. Similarly, the Backward Bifurcation Analyses offers insights into the diverse trajectories infectious diseases can undertake. Each of the section of this chapter is nuanced shift revealing the diverse patterns inherent in mathematical choreography. Sensitivity and Numerical Analysis shifts the focus to Sensitivity and Numerical Analysis of SEIR model. Sensitivity Analyses contribute to synthesizing art and science by translating infectious disease mysteries into mathematics. The chapter goes with Numerical Simulations, completing the intricate dance of mathematical precision and biological insights. Finally, this article, across these sections in it synthesizes art and science along with translating mathematics into a vibrant ballet to unravel infectious disease mysteries.

# Formulation and Analysis:

## 2.1: Overview

In this chapter we will formulate the model using classes of population and respective parameters. The well-posedness of the system of the model will be examined since it ensures that the system has a unique solution that is stable and continuous with respect to small changes in the initial conditions or parameters. The basic reproduction number will be examined with the help of next-generation matrix because it is important for predicting the epidemiology disease, evaluating the effectiveness of control measures, comparing with other diseases. After that we are going to analysis the stability of the model evaluating the equilibrium points. It is essential for assessing model validity and identifying critical thresholds.

## 2.2: Model Formulation

The dynamics of an infectious disease or any virus in modelling, we consider the entire population in which it occurs. The models we consider assume that N is the total size of population which is constant.
To formulate our model at each time t, we split the population N into four categories [6], [7], [8].

S: The susceptible class, or individuals who are not presently infected but could get the disease

E: The exposed class, a person who has been infected but not yet infectious.

I: The infective class, individuals who are infected with the disease and currently infectious.

R: The removed class, individuals who are immune to the disease by nature, or have fully recovered from it, or have passed away.

### 2.2.1: Equations for Modified SEIR Model

Model formulation is an important tool for mathematical modelling [6], [9]. To formulate models, we need to construct equations. Here, we will discuss equations for modified SEIR model [7]. To construct a model, we need to set a diagram to understand the assumption easily.

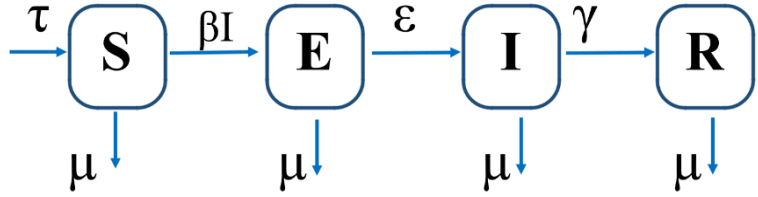

Figure 1: Schematic Diagram of SEIR Model

For considering the above situation, the amount of susceptible persons abates and hence is negative. The rate at which individuals enter exposed populations and depart the susceptible population is equal. The amount of individuals in exposed class E increases for the susceptible individuals to become exposed to the disease. The term ε will be the rate of being infectious from the exposed class. $\gamma$ be the rate at which an infected individual may be recovered and μR be the death rate of recovered individual.

**The corresponding model equations are as follows:**

$$\frac{dS}{dt} = \tau - \mu S - \beta SI$$
$$\frac{dE}{dt} = \beta SI - (\mu + \varepsilon)E$$
$$\frac{dI}{dt} = \varepsilon E - (\mu + \gamma)I \quad \quad (2.1)$$
$$\frac{dR}{dt} = \gamma I - \mu R$$

with initial conditions S(0) > 0, E(0) ≥ 0, I(0) > 0 and R(0) ≥ 0, whereas total population is as follows,

$$N(t) = S(t) + E(t) + I(t) + R(t) \quad \quad (2.2)$$

## 2.2.2: Parameter Identity

| Parameter | Interpretation |
|---|---|
| τ | Rate of Recruitment at which new individuals transfer to the susceptible class |
| μ | Rate of natural death |
| β | infection rate |
| ε | Rate of infectiousness of affected individuals |
| γ | Recovery rate |

Table 1: Parameter description

## 2.3: Positivity and boundedness of the solutions of the model

Here, under non-negative primary conditions, We will provide the positive arrangement of the system (2.1)
$N_0 = N(0)$, $S_0 = S(0)$, $E_0 = E(0)$, $I_0 = I(0)$, $R_0 = R(0)$,
are non-negative at any time $t$, let the solutions to the model are
$N = N(t)$, $S = S(t)$, $E = E(t)$, $I = I(t)$, $R = R(t)$,
Which are also non-negative i.e., the closed *region* $\Omega = \{(S, E, I, R) \in \mathbb{R}^4 +: S(0) > 0, E(0) \geq 0, I(0) > 0, R(0) \geq 0\}$. The boundedness and non-negative criteria of the solution confirm that the solution of the model exist, and the model has a unique solution [10].

From equation (2.2), we get, $N(t) = S(t) + I(t) + R(t) + E(t)$
Thence,

$$\frac{dN(t)}{dt} = \frac{dS(t)}{dt} + \frac{dE(t)}{dt} + \frac{dI(t)}{dt} + \frac{dR(t)}{dt} \quad \text{or,} \frac{dN(t)}{dt} = \tau - \mu S - \mu E - \mu I - \mu R$$

$$\text{or,} \frac{dN(t)}{dt} = \tau - \mu(S + E + I + R) \quad \text{or,} \frac{dN(t)}{dt} = \tau - \mu N$$

$$\text{or,} \frac{dN(t)}{dt} + \mu N = \tau \tag{2.3}$$

The different solutions of exposed, infected, and recovered class remain positive [9], [11], [12].

$$\frac{dE}{dt} = \beta SI \geq 0; \quad \frac{dI}{dt} = \varepsilon E \geq 0; \quad \frac{dR}{dt} = \gamma I \geq 0$$

From (2.3)

$$\frac{dN}{dt} = \tau - \mu N \leq \tau$$

It is evident that if $N > 0$ and $\frac{dN}{dt} < 0$,
Therefore, all the solutions of system (2.1) which have non-negative introduction qualities on the space $\mathbb{R}^4$, are bounded and exist within the interval $[0, \infty]$.

## 2.4: Calculating Basic Reproduction Number

The basic reproduction number is important for the analysis of infectious diseases. It is a threshold value that establishes whether the disease will keep happening or go away. It is described as the expected number of secondary cases that a typical infectious individual will cause. In this section, we used the next generation matrix to calculate the basic reproduction number[11], [12], [13], [14].
From the model equation,

$$\frac{dE}{dt} = \beta SI - (\mu + \varepsilon)E$$
$$\frac{dI}{dt} = \varepsilon E - (\mu + R)I$$

The Jacobean matrix for the system of the two infected equation is as follows:

$$J = \begin{bmatrix} -(\mu + \varepsilon) & \beta S \\ \varepsilon & -(\mu + R) \end{bmatrix}$$

The basic reproduction number is calculated from the largest positive eigenvalue of the next generation matrix $FV^{-1}$

where, $F = \begin{bmatrix} 0 & \beta S \\ \varepsilon & 0 \end{bmatrix}$ is Transmission matrix, $V = \begin{bmatrix} \mu + \varepsilon & 0 \\ 0 & \mu + \gamma \end{bmatrix}$ is Transition matrix

And the next generation matrix $FV^{-1}$ for system is,

$FV^{-1} = \begin{bmatrix} 0 & \frac{\beta S}{\mu + \gamma} \\ \frac{\varepsilon}{\mu + \gamma} & 0 \end{bmatrix}$. The maximum of eigenvalues of $FV^{-1}$ is $\frac{\varepsilon \beta}{(\mu+\gamma)(\mu+\varepsilon)}$. Therefore, the basic reproduction number, $R_0 = \frac{\varepsilon \beta}{(\mu+\gamma)(\mu+\varepsilon)}$ which is positive since all parameters are positive.

With the help of basic reproduction number, the steady state of the system can be discussed. The three situations of the steady state are as follows:

$R_0 < 1$ means a single primary infection causes less than one secondary infection. This implies the DFE is stable and the disease cannot persist in the population.

$R_0 > 1$ means a single primary infection causes more than one secondary infection. This implies the DFE is unstable. As a result, an epidemic breaks out.

$R_0 = 1$ means a single primary infection causes one secondary infection. This implies a trans-critical stage at disease free equilibrium (DFE).

## 2.5: Local Stability Analysis

Stability Analysis is an essential component of dynamical systems theory that explores the evolution of a system over time. Its goal is to ascertain whether the system's trajectories will converge to a specific equilibrium point, oscillate around it, or display chaotic behavior. Understanding the stability of a dynamical system is critical for predicting its long-term behavior and its response to various inputs or disturbances.

Throughout this section, we set up the stability by utilizing the basic reproduction number, provided the equilibrium point exists. This section also interprets the nature of the disease-free
equilibrium point and the nature of the EE point. To determine the equilibrium point, the steady-state condition is applied to the system of model equations (2.1)

At steady state, $\frac{dk}{dt} = 0$. where, $k \in \{S, E, I, R\}$. Therefore,

$$\frac{dS}{dt} = \tau - \mu S - \beta SI = 0$$
$$\frac{dE}{dt} = \beta SI - (\mu + \varepsilon)E = 0$$
$$\frac{dI}{dt} = \varepsilon E - (\mu + \gamma)I = 0 \qquad (2.4)$$
$$\frac{dR}{dt} = \gamma I - \mu R = 0$$

### 2.5.1: Stability of DFE (DFE)

When a disease stays for a long time in a region i.e. at t tends to infinite, it should be die out. The DFE requires that the infection should be zero i.e. $I(t) = 0$. From the last equation of system (2.4),

$$\frac{dR}{dt} = 0 \Rightarrow \gamma I - \mu R = 0 \Rightarrow I = \frac{\mu R}{\gamma}$$

For the DFE $I(t) = 0$. Therefore, $R = 0$

$$\frac{dE}{dt} = 0 \Rightarrow \beta SI - (\mu + \varepsilon)E = 0 \Rightarrow E = 0 \text{ . Also,}$$

$$\frac{dS}{dt} = 0 \Rightarrow \tau - \mu S - \beta SI = 0 \Rightarrow \tau - \mu S = 0 \Rightarrow S = \frac{\tau}{\mu}$$

Therefore, the set of the DFE point $E_o$ has been found from the solution set of the system of equation (2.4)

$E_{DFE} = (S_0, E_0, I_0, R_0) = \left(\frac{\tau}{\mu}, 0, 0, 0\right)$

Therefore, the local stability of DFE point have been estimated when the basic reproduction number $R_o > 1$ i.e. the disease will die out from the region.

**Theorem: 2.1**

If $R_o < 1$ then DFE stability point is locally asymptotically stable [9], [11], [12].

**Proof**: The DFE point is $E_{DFE} = (S_o, E_o, I_o, R_o) = \left(\frac{\tau}{\mu}, 0, 0, 0\right)$ and the basic reproduction number at that point is,

$R_o = \frac{\varepsilon \beta \tau}{(\mu+R)(\mu+\varepsilon)\mu}$

The Jacobian of the system of the DFE point is,

$$J_{DEF} = \begin{bmatrix} -\mu & 0 & \frac{-\beta\tau}{\mu} & 0 \\ 0 & -(\mu+\varepsilon) & \frac{\beta\tau}{\mu} & 0 \\ 0 & \varepsilon & -(\mu+\gamma) & 0 \\ 0 & 0 & \gamma & -\mu \end{bmatrix}$$

Characteristic equation of the Jacobian is,

$|J_{DFE} - \lambda I| = 0$, which gives,

$$\begin{vmatrix} -\mu-\lambda & 0 & \frac{-\beta\tau}{\mu} & 0 \\ 0 & -(\mu+\varepsilon)-\lambda & \frac{\beta\tau}{\mu} & 0 \\ 0 & \varepsilon & -(\mu+\gamma)-\lambda & 0 \\ 0 & 0 & \gamma & -\mu-\lambda \end{vmatrix} = 0$$

or, $(\mu+\lambda)^2 \left[\lambda^2 + (2\mu+\varepsilon+\gamma)\lambda + \left\{(\mu+\varepsilon)(\mu+\gamma) - \frac{\beta\tau\varepsilon}{\mu}\right\}\right] = 0$

or, $(\mu+\lambda)^2 (\lambda^2 + a_1\lambda + a_2) = 0$

Here, $a_1, a_2$ are the co-efficient of the characteristic's equation. For stability, they must be positive according to Descartes's rule of sign ("The number of positive real roots of a polynomial is equal to the number of variations in sign of its coefficients or less than that by an even integer. The number of negative real roots is similarly determined") [15]

$a_1 = 2\mu + \varepsilon + \gamma$

$a_2 = (\mu+\varepsilon)(\mu+\gamma) - \frac{\beta\tau\varepsilon}{\mu} = (\mu+\varepsilon)(\mu+\gamma)\left[1 - \frac{\beta\varepsilon\tau}{\mu(\mu+\varepsilon)(\mu+\gamma)}\right] = (\mu+\varepsilon)(\mu+\gamma)[1 - R_o]$

Using Descartes's rule of sign, all the roots are negative for $R_o > 1$ i.e. the system is locally asymptotically stable for $R_o < 1$ and it is unstable for $R_o > 1$. Hence proved.

2.5.2: Stability at the EE (EE)

We now prove the local stability of the (EE) point when $R_o > 1$. In a biological context, the disease will be endured in the population in case the basic reproduction number exceeds unity.

The determination of the EE state arises from the solution set of the system of equations (2.4)
When $I^* \neq 0$

$$\frac{dI^*}{dt} = 0 \Rightarrow \varepsilon E^* - (\mu + \gamma)I^* = 0 \Rightarrow E^* = \frac{(\mu + \gamma)I^*}{\varepsilon} \tag{2.5}$$

Now,

$$\frac{dE^*}{dt} = 0 \Rightarrow \beta S^* I^* - (\mu + \varepsilon)E^* = 0 \Rightarrow E^* = \frac{\beta S^* I^*}{\mu + \varepsilon} \tag{2.6}$$

So, from (2.5) and (2.6) we can write,

$$\therefore \frac{(\mu + \gamma)I^*}{\varepsilon} = \frac{\beta S^* I^*}{\mu + \varepsilon} \Rightarrow I^* \left\{ \frac{\mu + \gamma}{\varepsilon} - \frac{\beta S^*}{\mu + \varepsilon} \right\} = 0$$

Since $I^* \neq 0$, so we can write,

$$\frac{\mu + \gamma}{\varepsilon} - \frac{\beta S^*}{\mu + \varepsilon} = 0 \Rightarrow \frac{\beta S^*}{\mu + \varepsilon} = \frac{\mu + \gamma}{\varepsilon} \Rightarrow S^* = \frac{(\mu + \gamma)(\mu + \varepsilon)}{\beta \varepsilon} = \frac{1}{R_0}$$

So,

$$\frac{dS^*}{dt} = 0 \Rightarrow \tau - \beta S^* I^* - \mu S^* = 0 \Rightarrow \tau - \beta \cdot \frac{(\mu + \gamma)(\mu + \varepsilon)}{\beta \varepsilon} I^* - \mu \cdot \frac{(\mu + \gamma)(\mu + \varepsilon)}{\beta \varepsilon} = 0$$

$$\Rightarrow \frac{(\mu + \gamma)(\mu + \varepsilon)}{\varepsilon} I^* = \frac{\tau \beta \varepsilon - \mu(\mu + \gamma)(\mu + \varepsilon)}{\beta \varepsilon} \Rightarrow I^* = \frac{\tau \beta \varepsilon - \mu(\mu + \gamma)(\mu + \varepsilon)}{\beta(\mu + \gamma)(\mu + \varepsilon)}$$

$$= \frac{(\mu + \gamma)(\mu + \varepsilon)\left\{\frac{\gamma \beta \varepsilon}{(\mu + \gamma)(\mu + \varepsilon)} - \mu\right\}}{\beta(\mu + \gamma)(\mu + \varepsilon)} = \frac{1}{\beta}\{\tau R_0 - \mu\}$$

From (2.5),

$$E^* = \frac{(\mu + \gamma)}{\varepsilon} \cdot \frac{\tau \beta \varepsilon - \mu(\mu + \gamma)(\mu + \varepsilon)}{\beta(\mu + \gamma)(\mu + \varepsilon)} = \frac{(\mu + \gamma)(\mu + \varepsilon)\left[\mu - \frac{\tau \beta \varepsilon}{(\mu + \gamma)(\mu + \varepsilon)}\right]}{\beta \varepsilon(\mu + \varepsilon)}$$

$$= \frac{(\mu + \gamma)}{\beta \varepsilon}[\mu - \tau R_0] = \frac{\mu - \tau R_0}{R_0(\mu + \varepsilon)}$$

Then,

$$\frac{dR^*}{dt} = 0 \Rightarrow \gamma I^* - \mu R^* = 0 \Rightarrow R^* = \frac{\gamma}{\mu} I^* \Rightarrow R^* = \frac{\gamma[\tau \beta \varepsilon - \mu(\mu + \gamma)(\mu + \varepsilon)]}{\beta \mu(\mu + \gamma)(\mu + \varepsilon)}$$

$$= \frac{\gamma \tau \beta \varepsilon - \gamma \mu(\mu + \gamma)(\mu + \varepsilon)}{\beta \mu(\mu + \gamma)(\mu + \varepsilon)} = \frac{(\mu + \gamma)(\mu + \varepsilon)\left\{\frac{\gamma \tau \beta \varepsilon}{(\mu + \gamma)(\mu + \varepsilon)} - \gamma \mu\right\}}{\beta \mu (\mu + \gamma)(\mu + \varepsilon)}$$

$$= \frac{\gamma}{\beta \mu} \times \{\tau R_0 - \mu\}$$

$$\therefore (S^*, E^*, I^*, R^*)$$
$$= \left( \frac{(\mu + \gamma)(\mu + \varepsilon)}{\beta \varepsilon}, \frac{\tau \beta \varepsilon - \mu(\mu + \gamma)(\mu + \varepsilon)}{\beta \varepsilon(\mu + \varepsilon)}, \frac{\tau \beta \varepsilon - \mu(\mu + \gamma)(\mu + \varepsilon)}{\beta(\mu + \gamma)(\mu + \varepsilon)}, \frac{\gamma[\tau \beta \varepsilon - \mu(\mu + \gamma)(\mu + \varepsilon)]}{\beta \mu(\mu + \gamma)(\mu + \varepsilon)} \right)$$
$$= \left( \frac{1}{R_0}, \frac{\mu - \tau R_0}{R_0(\mu + \varepsilon)}, \frac{1}{\beta}\{\tau R_0 - \mu\}, \frac{\gamma}{\beta \mu} \times \{\tau R_0 - \mu\} \right)$$

is a critical point when $I^* \neq 0$. i.e. EE point.

**Theorem: 2.2**

The EE point $E_{EE}$ is locally asymptotically stable if $R_o > 1$ [9], [11], [12].

**Proof**: The Jacobian of the system at the EE point is,

$$J_{EE} = \begin{bmatrix} -\mu - \beta I^* & 0 & -\beta S^* & 0 \\ \beta I^* & -(\mu + \varepsilon) & \beta S^* & 0 \\ 0 & \varepsilon & -(\mu + \gamma) - \lambda & 0 \\ 0 & 0 & \gamma & -\mu \end{bmatrix} \text{ or } J_{EE} = \begin{bmatrix} a_1 & 0 & b_1 & 0 \\ a_2 & b_2 & c_2 & 0 \\ 0 & \varepsilon & c_3 & 0 \\ 0 & 0 & \gamma & -\mu \end{bmatrix}$$

The characteristics equation of the Jacobian is $|J_{EE} - \lambda I| = 0$, which gives,

$$\begin{vmatrix} a_1 - \lambda & 0 & b_1 & 0 \\ a_2 & b_2 - \lambda & c_2 & 0 \\ 0 & \varepsilon & c_3 - \lambda & 0 \\ 0 & 0 & \gamma & -\mu - \lambda \end{vmatrix} = 0$$

Therefore, the characteristic equation is, $(\lambda + \mu)(\lambda^3 + A\lambda^2 + B\lambda + C) = 0$

Where,

$A = -a_1 - b_2 - c_3 = -(-\mu - \beta I^*) + \mu + \varepsilon + \mu + \gamma = 3\mu + \varepsilon + \gamma + \beta I^* > 0$
$B = a_1 b_2 + a_1 c_3 + b_2 c_3 - c_2 \varepsilon = (\mu + \beta I^*)(2\mu + \varepsilon + \gamma) + (\mu + \varepsilon)(\mu + \gamma) - \beta S^* \varepsilon > 0$
$C = a_1 \varepsilon c_2 - b_1 \varepsilon a_2 - a_1 b_2 c_3 = -(\mu + \beta I^*) \varepsilon \beta S^* + \beta S^* \varepsilon \beta I^* + (\mu + \beta I^*)(\mu + \varepsilon)(\mu + \gamma)$
$\phantom{C} = (\mu + \beta I^*)[(\mu + \varepsilon)(\mu + \gamma) - \beta \varepsilon S^*] + \beta^2 S^* I^* \varepsilon > 0$

By Ruth-Hurwitz condition ("For a polynomial with real coefficients, the system described by that polynomial is stable if and only if all coefficients are positive, and the determinants of the Hurwitz matrices, formed from the coefficients, are all positive."), [15]
A. $B - C > 0$

## 2.6: Global Stability by Lyapunov Concept

To determine global stability, we use the Lyapunov function and the stability theorem [16]. Before showing the function and theorem we need to introduce some terms.

**Positive Definite function:** A function F: $\mathbb{R}^n \to \mathbb{R}$ is called positive definite function in a neighborhood U of origin if the functional value of F is zero at origin and the functional value of F is greater than zero for all values of X in U except origin.

**Negative Definite function:** A function F: $\mathbb{R}^n \to \mathbb{R}$ is called positive definite function in a neighborhood U of origin if the functional value of F is zero at origin and the functional value of F is less than zero for all values of X in U except origin.

**Positive Semi Definite function:** A function F: $\mathbb{R}^n \to \mathbb{R}$ is called positive definite function in a neighborhood U of origin if the functional value of F is zero at origin and the functional value of F is greater than or equal to zero for all values of X in U except origin.

**Negative Semi Definite function:** A function F: $\mathbb{R}^n \to \mathbb{R}$ is said to be positive definite function in a neighborhood U of origin if the functional value of F is zero at origin and the functional value of F is less than or equal to zero for all values of X in U except origin.

A scalar function defined on phase space that can be used to demonstrate the stability of an equilibrium point is called a **Lyapunov function** [17].

Assume that V(X) is a continuously differentiable function in the neighborhood U of the origin. The function V(X) is called the Lyapunov function for an autonomous system X' = f(x) if the following conditions hold

1. For every X in U\ {0}, V(X) >0, that is V(X) is positive definite.
   2. For all X in U, (dV/dt) < 0, that is dV/dt is negative definite.
   3. V (0) is equal to 0.

This concept comes as a theorem which we call Lyapunov stability theorem [16], [18].

## 2.6.1: Lyapunov stability theorem

Consider a system X' = f(x) which has a equilibrium point at origin (0, . . . , 0) And there exist a bounded neighborhood N of the equilibrium point (0, . . . , 0) and a function V is defined in N which is a Lyapunov function for this system then this theorem states that any point either the any closed region or outside the closed region of that equilibrium point converges to that equilibrium point as time tends to infinite. i.e. if the following conditions are satisfied, then the equilibrium point is globally stable.

1. The first order partial derivatives of the Lyapunov function V are continuous in the neighborhood N.
2. V is positive definite.
3. dV/dt is negative semi definite.

Where,
$$\frac{dV}{dt} = \frac{\partial V}{\partial x_1}\frac{dx_1}{dt} + \frac{\partial V}{\partial x_2}\frac{dx_2}{dt} + \cdots + \frac{\partial V}{\partial x_n}\frac{dx_n}{dt}$$

So, the procedure to check global stability of a system by Lyapunov sense, first we need to get the equilibrium points of the system and check whether the system has an equilibrium point at origin. If the system has the equilibrium point at origin, then we need to construct the Lyapunov function, V and then check whether the function satisfies the Lyapunov theorem conditions. If the function satisfied all the conditions of Lyapunov theorem then we can say our system is globally stable at the origin. Here we discussed the procedure of Lyapunov to check the global stability at origin if the system has equilibrium point at origin. To check the stability of other equilibrium points except origin by Lyapunov procedure we need to transform the equilibrium points to origin so that we can check the stability by this procedure.

**To better understand the Lyapunov function, let's examine an example:**

Using the Lyapunov function approach, different differential equations and systems are examined for stability. For instance, let's consider an autonomous system [17].

$X' = f(x)$, or $\frac{dx_i}{dt} = f_i(x_1, x_2, x_3, \ldots x_n)$, In this case, i = 1, 2,...,n

With x≡0, which represents zero equilibrium.

Let the function V(X) = V (x$_1$, x$_2$, …, x$_n$) be continuously differentiable which is given in the neighborhood U of the origin. Allowing V(X) to represent the sum of all X in U\{0} and V(0) at the origin. For instance, these are functions of the kind.

$V(x_1, x_2) = ax_1^2 + bx_2^2$ or $V(x_1, x_2) = ax_1^2 + bx_2^4$, where a and b are positive constants.

Total derivative of the function of V(X) with respect to time t can be determined as follows:
$$\frac{dV}{dt} = \frac{\partial V}{\partial x_1}\frac{dx_1}{dt} + \frac{\partial V}{\partial x_2}\frac{dx_2}{dt} + \cdots + \frac{\partial V}{\partial x_n}\frac{dx_n}{dt}$$

This equation can be expressed as follows, just as a scalar product of two vectors:
$$\frac{dV}{dt} = grad\ V \cdot \frac{dX}{dt}$$

Where, $grad\ \text{V} = \left(\frac{\partial V}{\partial x_1}, \frac{\partial V}{\partial x_2}, \ldots, \frac{\partial V}{\partial x_n}\right)$ , $\frac{dX}{dt} = \left(\frac{dx_1}{dt}, \frac{dx_2}{dt}, \ldots \frac{dx_n}{dt}\right)$

The first vector is always going in the direction of the greatest rise in V(X), since it is the gradient of V(X). When a function V(X) is given |X|→∞, it usually grows with the distance

from the origin. The velocity vector comes next in the scalar product. It is always tangent to the phase trajectory.

Examine the following case: V(X) has a negative derivative in the neighborhood U of the origin.

$$\frac{dV}{dt} = \left(grad\ V \cdot \frac{dX}{dt}\right) < 0$$

This indicates that the gradient and velocity vectors have an angle φ greater than 90 degrees. A phase trajectory is going to move toward the origin if the derivative dV/dt is always negative, indicating that the system is stable. The system is unstable whenever the derivative dV/dt is positive because the trajectory goes away from the origin.

**Checking global stability of DFE point of SEIR model by Lyapunov theorem:**

$$\frac{dS}{dt} = \tau - \mu S - \beta SI$$
$$\frac{dE}{dt} = \beta SI - (\mu + \varepsilon)E$$
$$\frac{dI}{dt} = \varepsilon E - (\mu + \gamma)I \quad (2.7)$$
$$\frac{dR}{dt} = \gamma I - \mu R$$

Where, $s > 0,\ E \geqslant 0,\ I \geqslant 0,\ R \geqslant 0$

Disease free equilibrium point is $(S^*, E^*, I^*, R^*) = \left(\frac{\tau}{\mu}, 0, 0, 0\right)$.

We calculated this before.

Now to check the global stability transform the equilibrium point to origin.

$$\text{i.e.}\ \left(S^* - \frac{\tau}{\mu}, 0, 0, 0\right) = (0, 0, 0, 0)$$

Now let,

$$x = S^* - \frac{\tau}{\mu},\quad y = E^*,\quad z = I^*,\quad w = R^* \Rightarrow S^* = x + \frac{\tau}{\mu}$$

Then the system transforms to

$$\frac{d\left(x + \frac{\gamma}{\mu}\right)}{dt} = \tau - \mu\left(x + \frac{\gamma}{\mu}\right) - \beta\left(x + \frac{\gamma}{\mu}\right)(z)$$
$$\frac{dy}{dt} = \beta\left(x + \frac{\gamma}{\mu}\right)z - (\mu + \varepsilon)y$$
$$\frac{dz}{dt} = \varepsilon y - (\mu + \gamma)z$$
$$\frac{d\omega}{dt} = \gamma z - \mu w$$

$$\Rightarrow \begin{aligned}\frac{dx}{dt} &= -\mu x - \beta xz - \frac{\beta\gamma}{\mu}z \\ \frac{dy}{dt} &= \beta rz + \frac{\beta y}{\mu}z - (\mu + \varepsilon)y \\ \frac{dz}{dt} &= \varepsilon y - (\mu + \gamma)z \\ \frac{dw}{dt} &= \gamma z - \mu w\end{aligned} \quad (2.8)$$

so, the system transformed to (2.8) from (2.7) and the DFE point of (2.7) transformed to (0,0,0,0)

Now, consider the **Lyapunov function,**
$$V(x, y, z, w) = x + y + z + w \tag{2.9}$$

Now,

**(i)** The first partial derivatives of (2.9) are continuous. Since,
$$\frac{\partial V}{\partial x} = 1, \frac{\partial V}{\partial y} = 1, \frac{\partial V}{\partial z} = 1, \frac{\partial V}{\partial w} = 1\ ;$$

**(ii)** Now, $x > 0, y \geqslant 0, z \geqslant 0, w \geqslant 0$ since $s > 0, E \geqslant 0, I \geqslant 0, R \geqslant 0$.

And $V(x, y, z, w) = 0$ at $(0,0,0,0)$ and $V(x, y, z, w) > 0$, for all $(x, y, z, w)$ except $(x, y, z, w) = (0,0,0,0)$ since, $x > 0, y \geqslant 0, z \geqslant 0, w \geqslant 0$. So, $V$ is positive definite.

**(iii)** The time derivative of $v$ is
$$\frac{dV}{dt} = \frac{\partial V}{\partial x} \cdot \frac{dx}{dt} + \frac{\partial V}{\partial y} \cdot \frac{dy}{dt} + \frac{\partial V}{\partial z} \cdot \frac{dz}{dt} + \frac{\partial V}{\partial w} \cdot \frac{dw}{dt} \Rightarrow \frac{dV}{dt} = 1 \cdot \frac{dx}{dt} + 1 \cdot \frac{dy}{dt} + 1 \cdot \frac{dz}{dt} + 1 \cdot \frac{dw}{dt}$$
$$\Rightarrow \frac{dV}{dt} = \frac{dx}{dt} + \frac{dy}{dt} + \frac{dz}{dt} + \frac{dw}{dt} \tag{2.10}$$

now substituting the values of $\frac{dx}{dt}, \frac{dy}{dt}, \frac{dz}{dt}, \frac{dw}{dt}$ in (2.10) from (2.8), we get,

$\frac{dV}{dt} = -\mu x - \mu y - \mu z - \mu w \Rightarrow \frac{dV}{dt} = -\mu(x + y + z + w)$

Now, $\frac{dV}{dt} = 0$ at $(x, y, z, w) = (0,0,0,0)$ and $\frac{dV}{dt} \leqslant 0$, since $x + y + z + w \geqslant 0$.

So, $\frac{dV}{dt}$ is negative semi definite.

So, all the three conditions of Lyapunov theorem are satisfied. so, the DFE points of (2.5), $\left(\frac{\tau}{\mu}, 0,0,0\right)$ is **globally stable**

# Bifurcation Analysis and its Classification:

## 3.1: Overview

In this section we will describe some simple examples of bifurcations that occur for nonlinear systems. It is often used to describe the behavior of dynamic systems, such as in chaos theory and nonlinear dynamics, where small changes in parameters can lead to significant and unpredictable changes in the system's behavior. Bifurcation can result in the emergence of new patterns, behaviors, or states in a system, and is a key concept in understanding complex and nonlinear systems. It has two types:

**Local bifurcations** is a phenomenon in dynamical systems theory where a qualitative change in the behavior of a system occurs as a parameter is varied, typically near a fixed point or equilibrium point.

**Global bifurcation** is a phenomenon in dynamical systems theory where a qualitative change in the behavior of a system occurs across a wide range of parameter values, affecting the system as a whole.

## 3.2: Bifurcation in one dimension

### Saddle-node bifurcation:

In dynamical systems theory, a saddle-node bifurcation is a type of bifurcation where a fixed point of a system loses stability and disappears as parameter is varied. This phenomenon happens when the critical equilibrium of an autonomous system contains one zero eigenvalue.[19], [20], [21]. Now, let's consider the dynamical system defined by,
$$x' = x^2 + a; \ a \text{ is real.}$$
$$y' = -y.$$

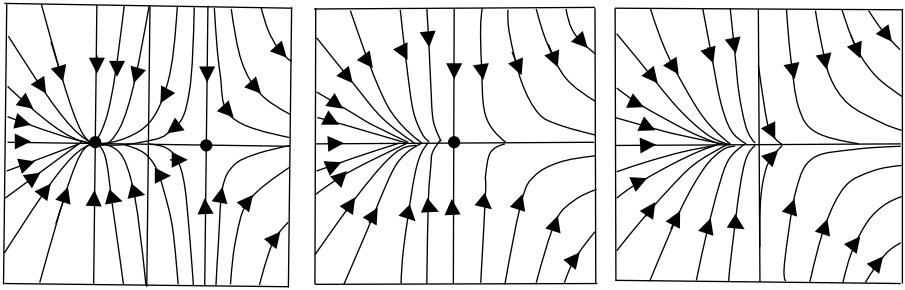

Figure 2: Bifurcation Diagram for a Saddle-node bifurcation [20]

To analyze the transcritical bifurcation in this system, we first find the fixed points by setting x' and y' equal to zero. For the x equation, we have $x^2 + a = 0$ which gives us two fixed points: $x = \sqrt{-a}$ and $= -\sqrt{-a}$. For the y equation, the fixed point is $y = 0$. Next, we analyze the stability of these fixed points by linearizing the system around each fixed point. We compute the Jacobian matrix:
$$J = \begin{bmatrix} 2x & 0 \\ 0 & -1 \end{bmatrix}$$
Evaluating the Jacobian at the fixed points, we find that the fixed point $(\sqrt{-a}, 0)$ is stable (both eigenvalues are negative) and the fixed point $(-\sqrt{-a}, 0)$ is unstable (one eigenvalue is positive). As the parameter $a$ crosses the critical value of zero, the fixed points at $(\sqrt{-a}, 0)$ and $(-\sqrt{-a}, 0)$ exchange stability. This exchange of stability characterizes the transcritical bifurcation in this system. For $a < 0$ trajectories near the stable fixed point will converge towards it, while trajectories near the unstable fixed point will move away from it.

At $a = 0$, the transcritical bifurcation occurs, causing the fixed points to exchange stability. The fixed point at $(\sqrt{0}, 0) = (0,0)$ becomes unstable, and the fixed point at $(-\sqrt{0}, 0) = (0,0)$ becomes stable. Trajectories near the origin will exhibit different behaviour before and after the bifurcation.

For $a > 0$ after the bifurcation, the fixed point at $(-\sqrt{-a}, 0)$ becomes stable, and the fixed point at $(\sqrt{-a}, 0)$ becomes unstable. Trajectories near the stable fixed point will converge towards it, while trajectories near the unstable fixed point will move away from it.

**Transcritical bifurcation**

A transcritical bifurcation is a type of bifurcation in dynamic systems where two families of fixed points exchange stability as a parameter is varied [19]. This leads to a qualitative change in the system's behavior, affecting stability and the types of solutions that can exist.

Considering the following dynamical system,
$$\frac{dx}{dt} = ax - bx^2, \text{ for } x, a, b \text{ real.}$$

To analyse this system, we first find the fixed points by setting $\frac{dx}{dt} = 0 \Rightarrow ax - bx^2 = 0 \Rightarrow x(a - bx) = 0$ The two fixed points are $x = 0$ and $= \frac{a}{b}$.

The Jacobian matrix of the system, $J = \left[\frac{d(ax-bx^2)}{dx}\right] = a - 2bx$

Evaluating the Jacobian matrix at the fixed points, At $x = 0, J(0) = a$ At $x = \frac{a}{b}, J\left(\frac{a}{b}\right) = -a$. The eigenvalues of the Jacobian matrix at $x = 0$ are both positive if $a > 0$, indicating an unstable fixed point. The eigenvalues at $x = \frac{a}{b}$ are both negative if $a > 0$, indicating a stable fixed point. Therefore, as the parameter $a$ is varied, the fixed points $x = 0$ and $x = \frac{a}{b}$ exchange stability at $a = 0$, which is the transcritical bifurcation point in this system.

**The Pitchfork bifurcation**

A pitchfork bifurcation is a type of bifurcation in dynamical systems where, as a parameter varied, a stable equilibrium point undergoes a symmetry-breaking transition, resulting in the emergence of two new equilibrium points – one stable and one unstable. The bifurcation is named after its characteristic shape resembling a pitchfork [19], [20].

Consider the dynamical system,
$$x' = x^3 - ax$$

There are three equilibria for this equation, at $x = 0$ and $x = \pm\sqrt{a}$ when $a > 0$. When $a \leq 0$, $x = 0$ is the only equilibrium point.

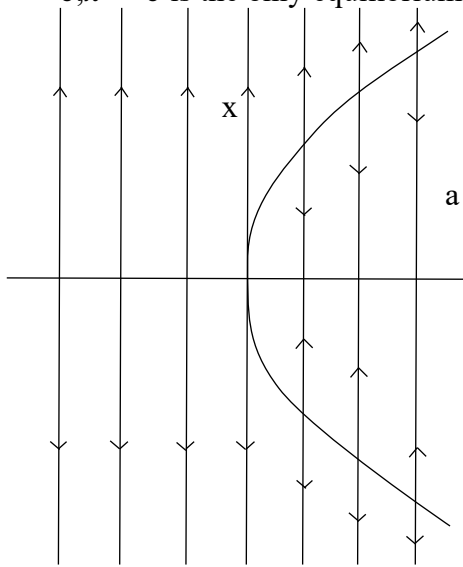

Figure 3: Pitchfork Bifurcation Diagram [20]

## 3.3: Bifurcation in two dimensions

**Hopf bifurcation**

**Definition**: A Hopf bifurcation is a type of bifurcation in dynamic systems where a stable equilibrium point loses stability, and a limit cycle is born as a parameter is varied [19]. This bifurcation is associated with the transition from a stable fixed point to a periodic oscillatory behavior in the system. It occurs when a pair of complex conjugate eigenvalues of the Jacobian matrix crosses the imaginary axis.

**Hopf bifurcation theorem**

Consider the two-dimensional system,
$$\begin{aligned} \frac{dx}{dt} &= f(x, y, \tau) \\ \frac{dy}{dt} &= g(x, y, \tau) \end{aligned} \quad (3.1)$$

where $\tau$ is the parameter and suppose that $(x(\tau), y(\tau))$ is the equilibrium point and $\alpha(\tau) \pm i\beta(\tau)$ are the eigenvalues of the Jacobian matrix which is evaluated at the equilibrium point. The equilibrium lies at the origin after the system is first transformed and the parameter $\tau$ at $\tau^* = 0$ gives purely imaginary eigenvalues. System (2.1) is rewritten as follows:

$$\frac{dx}{dt} = a_{11}(\tau)x + a_{12}(\tau)y + f_1(x, y, \tau)$$
$$\frac{dy}{dt} = a_{21}(\tau)x + a_{22}(\tau)y + g_1(x, y, \tau)$$
(3.2)

The linearization of the system (2.10) about the origin is given by
$\frac{dX}{dt} = J(\tau)X$, where $X = \frac{x}{y}$ and

$$J(\tau) = \begin{bmatrix} a_{11}(\tau) & a_{12}(\tau) \\ a_{21}(\tau) & a_{22}(\tau) \end{bmatrix}$$ is the matrix of Jacobian evaluated at the origin?

Example: Consider the two-dimensional system
$$x' = ax - y - x(x^2 + y^2)$$
$$y' = x + ay - y(x^2 + y^2).$$

The system contains an equilibrium point at the origin and the linearized form of the above system is
$$X' = \begin{pmatrix} a & -1 \\ 1 & a \end{pmatrix} X.$$

Computing the eigenvalues, $\lambda_1 = a + i$, and $\lambda_2 = a - i$ which are complex conjugate and crosses the imaginary axis. Therefore, solution of the system bifurcates at the origin. Transforming the system in polar co-ordinates by putting $x = r\cos\theta$, $y = r\sin\theta$ and

$$x' = \frac{dx}{dt} = \frac{dx}{dr} \cdot \frac{dr}{dt} + \frac{dx}{d\theta} \cdot \frac{d\theta}{dt} = r'\cos\theta - r\sin\theta \cdot \theta'$$
$$y' = \frac{dy}{dt} = \frac{dy}{dr} \cdot \frac{dr}{dt} + \frac{dy}{d\theta} \cdot \frac{d\theta}{dt} = r'\sin\theta + r\cos\theta \cdot \theta'$$

Now,
$$r'\cos\theta - r\sin\theta \cdot \theta' = ar\cos\theta - r\sin\theta - r^3\cos\theta \quad (3.3)$$
$$r'\sin\theta + r\cos\theta \cdot \theta' = ar\sin\theta + r\cos\theta - r^3\sin\theta \quad (3.4)$$

Solving (3.3) and (3.4) we get,
$$r' = ar - r^3$$
$$\theta' = 1.$$

Since $\theta'$ is non-zero the only equilibrium point of the system is the origin. For $a < 0$, $r' < 0$ for all $r > 0$ which represent the origin is a sink. For $a > 0$, the origin becomes a source. When $a > 0$, $r' = 0$ shows a limit cycle of radius $\sqrt{a}$ and a periodic solution with period $2\pi$. We also have $r' > 0$ if $0 < r < \sqrt{a}$, while $r' < 0$ if $r > \sqrt{a}$. Thus all trajectories spiral toward this circular solution as $t \to \infty$.

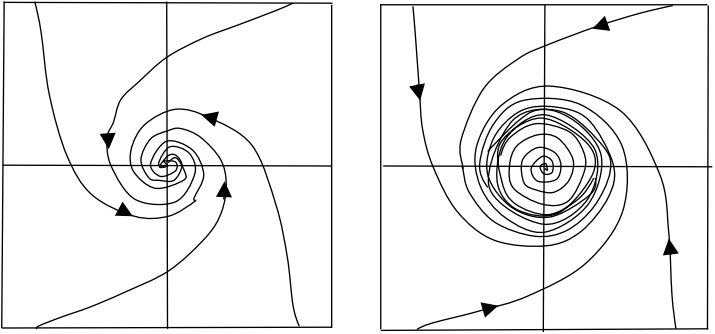

Figure 4: Hopf bifurcation for a<0 and a>0.

This type of bifurcation is called a *Hopf bifurcation.*
## 3.4: Forward Bifurcation Analysis

In this part, we will examine the type of the bifurcation involving the DFE $(S^*, E*, I^*)$ for $R_0 = 1$. The condition imposed on parameter values causes forward bifurcation to occur. Consider a general system of ODES with a parameter $\phi$;

$$\frac{dx}{dt} = f(x, \phi); f: \mathbb{R}^n \times \mathbb{R}^n, f \in c^2(\mathbb{R}^n \times \mathbb{R}^n) \tag{3.5}$$

It is assumed that 0the origin s an equilibrium for system (3.5) for all values of parameter $\phi$, that is, $f(0, \phi) \equiv 0$, for all $\phi = 0$
**Lemma-3.1** [22]: *(Castillo-Chavez and Song) Assume*
*(A1) $Q = D_x f(0,0)$ is the linearization matrix of system (3.5) around the equilibrium $x = 0$ with $\phi$ evaluated at 0. Zero is a simple eigenvalue of Q, and all other eigenvalues of Q have negative real parts.*
*(A2) Matrix Q has a (non-negative) right eigenvector w and a left eigenvector v corresponding to the zero eigenvalue.*
Let $f_k$ denotes the k-th component of f and,

$$\boldsymbol{a} = \sum_{k,i,j=1}^{3} v_k w_i w_j \frac{\partial^2 f_k}{\partial x_i \partial x_j}(0,0), \qquad \boldsymbol{b} = \sum_{k,i=1}^{3} v_k w_i \frac{\partial^2 f_k}{\partial x_i \partial \beta^*}(0,0).$$

Introducing $= x_1, E = x_2, I = x_3$, the system (2.1) becomes,

$$\frac{dx_1(t)}{dt} = \tau - \beta x_1(t) x_3(t) - \mu x_1(t) = f_1$$
$$\frac{dx_2(t)}{dt} = \beta x_1(t) x_3(t) - (\mu + \varepsilon) x_2(t) = f_2 \tag{3.6}$$
$$\frac{dx_3(t)}{dt} = \varepsilon x_2(t) - \mu x_3(t) - \gamma x_3(t) = f_3$$

We will apply Lemma (3.1) to show that system (3.6) may exhibit backward bifurcation when $R_0 = 1$. Taking $\tau = \mu$, we consider the DFE $\left(\frac{\tau}{\mu}, 0, 0\right)$ or $(1,0,0)$ and observe that the condition $R_0 = 1$ can be seen in terms of the parameter as $\beta = \beta^* = \frac{(\mu+\varepsilon)(\mu+\gamma)}{\varepsilon}$
The eigenvalues of the matrix,

$$J_{at\beta=\beta^*} = \begin{pmatrix} -\mu & 0 & -\beta \\ 0 & -\mu-\varepsilon & \beta \\ 0 & \varepsilon & -(\mu+\gamma) \end{pmatrix}$$

are given by, $\lambda_1 = -\mu, \lambda_2 = -(2\mu + \gamma + \varepsilon), \lambda_3 = 0$.
Evaluating the Jacobian matrix, we found three eigenvalues, one is zero and the other two eigenvalues are real and negative. Hence, when $\beta = \beta^*$ (or equivalently when $R_0 = 1$), the DFE is a non-hyperbolic equilibrium: the assumption (1) of Lemma is then verified.
Denoting $\boldsymbol{w} = (w_1, w_2, w_3)^\top$ a right eigenvector associated with the zero eigenvalue $\lambda_3 = 0$. It is found that,

$$\begin{pmatrix} -\mu & 0 & -\beta \\ 0 & -\mu-\varepsilon & \beta \\ 0 & \varepsilon & -(\mu+\gamma) \end{pmatrix} \begin{pmatrix} w_1 \\ w_2 \\ w_3 \end{pmatrix} = 0$$

Thus, we can get,

$$\begin{cases} -\mu w_1 - \beta w_3 = 0 \\ (-\mu - \varepsilon)w_2 + \beta w_3 = 0 \\ \varepsilon w_2 - (\mu + \gamma)w_3 = 0 \end{cases}$$

This implies $w_1 = -\left(\frac{\beta}{\mu}\right)w_3, w_2 = \frac{\beta}{(\mu+\varepsilon)}w_3, w_3 = w_3 > 0$

Therefore, the right eigenvector is

$$\mathbf{w} = \left(-\frac{\beta}{\mu}, \frac{\beta}{(\mu+\varepsilon)}, 1\right)$$

Furthermore, the left eigenvector $\mathbf{v} = (v_1, v_2, v_3)$ satisfying $\mathbf{v}.\mathbf{w} = 1$ is given by

$$-\mu v_1 = 0$$
$$-(\mu + \varepsilon)v_2 + \varepsilon v_3 = 0$$
$$-\frac{\beta}{\mu}v_1 + \frac{\beta}{\mu}v_2 - (\mu + \gamma)v_3 = 0$$

Then, the left eigenvector $v$ becomes,

$$\mathbf{v} = \left(0, \frac{(\mu+\gamma)}{\beta}, 1\right)$$

Evaluating the partial derivatives at the DFE, we obtain,

$$\frac{\partial^2 f_1}{\partial x_1 \partial x_3} = \frac{\partial^2 f_1}{\partial x_3 \partial x_1} = -\beta$$

$$\frac{\partial^2 f_2}{\partial x_1 \partial x_3} = \frac{\partial^2 f}{\partial x_3 \partial x_1} = \beta$$

Thus, we can compute the coefficients $\mathbf{a}$ and $\mathbf{b}$ defined in **Lemma 3.1**, i.e.

$$\mathbf{a} = \sum_{k,i,j=1}^{3} v_k w_i w_j \frac{\partial^2 f_k(0,0)}{\partial x_i \partial x_j}$$

$$\mathbf{b} = \sum_{k,i=1}^{3} v_k w_i \frac{\partial^2 f_k(0,0)}{\partial x_i \partial \beta^*}$$

Considering system (3.6) and the co-efficient in $\mathbf{a}$ and $\mathbf{b}$ only the nonzero derivatives of the terms $\frac{\partial^2 f_k(0,0)}{\partial x_i \partial x_j}$ and $\frac{\partial^2 f_k(0,0)}{\partial x_i \partial \beta^*}$, it follows that

$$\mathbf{a} = \sum_{k,i,j}^{3} v_k w_i w_j \frac{\partial^2 f_k(0,0)}{\partial x_i \partial x_j} = v_2 w_1 w_3 \beta > 0 \qquad (3.7)$$

Choosing $\beta^*$ as the bifurcation parameter, and we evaluate associated partial derivatives of $f_i$ are $\frac{\partial^2 f_1}{\partial x_3 \partial \beta^*} = -1;$

$\frac{\partial^2 f_2}{\partial x_3 \partial \beta^*} = 1,$ so we get,

$$\mathbf{b} = \sum_{k,i=1}^{3} v_k w_i \frac{\partial^2 f_k(0,0)}{\partial x_i \partial \beta^*}(f_i) = \frac{\mu+\gamma}{\beta}w_3 > 0 \qquad (3.8)$$

The following result can be determined based on the coefficients of $\mathbf{a}$ and $\mathbf{b}$.

**Theorem 3.1** [22]: *The system (3.6) exhibits backward bifurcation at $R_0 = 1$ whenever the bifurcation coefficients denoted by $\mathbf{a}$ and $\mathbf{b}$ given equations (3.7) and (3.8) are positive. However, if the co-efficient $\mathbf{a}$ is negative from equation (a), then the system (3.6) will not undergo a backward bifurcation at $R_0 = 1$.*

## 3.5: Backward Bifurcation Analysis

One of the most important requirements for models of infection control is having the connected $R_0$. In this part, the process of backward bifurcation for the following system is examined [12]. The system is expressed as $\frac{dx}{dt} = f(x)$, with $f(x) = (f_1(x), f_2(x), f_3(x), f_4(x))$ and

$$\frac{dx_1}{dt} = f_1 = -(\beta_1 + \beta_2)x_1 x_3 + \varepsilon x_2 + \alpha x_4$$
$$\frac{dx_2}{dt} = f_2 = \beta_1 x_1 x_3 - \varepsilon x_2 - \phi x_2 - \sigma x_2$$
$$\frac{dx_3}{dt} = f_3 = \beta_2 x_1 x_3 - \gamma x_3 - \delta x_3 + \phi x_2$$
$$\frac{dx_4}{dt} = f_4 = \gamma x_3 - \alpha x_4$$

(3.9)

Let $x = (x_1, x_2, x_3, x_4)^\top$. We use the center manifold theorem[23], to do a bifurcation examination. In terms of the bifurcation parameter $\beta_2$, when $R_0 = 1$, then $\beta_2 = \beta_2^* = \frac{(\gamma+\delta)(\varepsilon+\phi+\sigma)-\beta_1\phi}{\varepsilon+\phi+\sigma}$

At the DFE $E_0$ evaluated for $\beta_2 = \beta_2^*$, the system (3.7) has a simple eigenvalue with zero real components, while all other eigenvalues have a negative real part. Therefore, we examine the dynamics of (3.7) at $\beta_2 = \beta_2^*$ using the centre manifold theorem. The Jacobian of (3.7) at $\beta_2 = \beta_2^*$ is represented by $[J_{E_0}]_{\beta_2=\beta_2^*}$ and contains an eigenvector on the right (equivalent to the zero eigenvalues) given by $w = (w_1, w_2, w_3, w_4, \omega_5)^\top$.

In this case, $w_1 = w_1 > 0, w_2 = \frac{\beta_1 s^0}{\varepsilon+\phi+\sigma}, w_3 = 1, w_4 = \frac{\gamma}{\alpha}, w_5 = w_5 > 0$. Similar to this, we get another eigenvector on the left, $v = (v_1, v_2, v_3, v_4, v_5)^\top$ (corresponding to the zero eigenvalues) from $[J_{E_0}]_{\beta_2=\beta_2^*}$, $v_1 = 0; v_2 = 1; v_3 = \frac{\varphi}{\varepsilon+\varphi+\sigma}, v_4 = 0; v_5 = 0$. We compute the second-order partial derivatives of $f_i$ at the disease-free equilibrium in the following way to show that a backward bifurcation exists. The second-order partial derivatives of $f_i$ at the DFE as follows:

$$\frac{\partial^2 f_1}{\partial x_1 \partial x_3} = -(\beta_1 + \beta_2); \frac{\partial^2 f_2}{\partial x_1 \partial x_3} = \beta_1; \frac{\partial^2 f_3}{\partial x_1 \partial x_3} = \beta_2;$$

$$\frac{\partial^2 f_1}{\partial x_1 \partial \beta_2} = -x_3; \frac{\partial^2 f_1}{\partial x_3 \partial \beta_2} = -x_1; \frac{\partial^2 f_3}{\partial x_1 \partial \beta_2} = x_3; \frac{\partial^2 f_3}{\partial f_3 \partial f_2} = x_1$$

Now, we find the coefficients of $a$ and $b$ defined as follows:

$$a = \sum_{k,i,j=1}^{4} v_k w_i w_j \cdot \frac{\partial^2 f_k(0, \beta_2^*)}{\partial x_i \partial x_j} = w_1 \left[\frac{\varphi \beta_1}{\varepsilon + \phi + \sigma} + \beta_2\right]$$

$$\text{and } b = \sum_{k,i=1}^{5} v_k w_i \frac{\partial^2 f_k(0,0)}{\partial x_i \partial \beta_2} = \frac{\varepsilon + \varphi + \sigma}{\varphi} > 0$$

system (3.9) undergoes a backward bifurcation at $R_0 = 1$ because the coefficient $b$ is always positive [12].

# Sensitivity Analysis and Numerical Simulations:

In this section we like to address the issues that is key to spread the diseases in the society based on the existing parameters of the system. The we visualize the outcome of the system as solution trajectories.

## 4.1: Sensitivity Analysis

Sensitivity analysis helps the researcher to determine the effect of the parameter in the model. It is a technique to evaluate how the result will change based on changing input variables [24-30]. Here, we study the influence on basic reproduction number by the variation of the parameters in the expression of basic reproduction number. Denoting basic reproduction number by $R_0$, the sensitivity can be measured by the partial derivative with respect to the parameters of $R_0$ and it is represented by the notation, $S_P = \frac{\partial R_0}{\partial P} \frac{P}{R_0}$
where $p$ is any parameter in the expression of $R_0$.

To analyze sensitivity, we first consider the system (2.1)

$$\frac{dS}{dt} = \tau - \mu S - \beta SI$$

$$\frac{dE}{dt} = \beta SI - (\mu + \varepsilon)E$$

$$\frac{dI}{dt} = \varepsilon E - (\mu + \gamma)I$$

$$\frac{dR}{dt} = \gamma I - \mu R$$

From this system we get the basic reproduction number $R_0 = \frac{\varepsilon \beta}{(\mu+\gamma)(\mu+\varepsilon)}$

Where,

| Parameter | Identity |
|---|---|
| $\beta$ | Infection Rate |
| $\varepsilon$ | Rate at which infected become infectious |
| $\gamma$ | Recovery rate |
| $\mu$ | Natural death rate |

Table 2

Now we assume the parameter value as following:

| Parameter | Value | Data source |
|---|---|---|
| $\beta$ | 0.25 | assumed |
| $\varepsilon$ | 0.06 | assumed |
| $\gamma$ | 0.07 | assumed |
| $\mu$ | 0.005 | assumed |

Table 3

We see the parameters appearing in the expression of basic reproduction number are $\beta, \varepsilon, \mu$ and $\gamma$.

Starting with $\beta$,

$$S_\beta = \frac{\partial R_0}{\partial \beta}\frac{\beta}{R_0} = \frac{\partial}{2\beta}\frac{\varepsilon\beta}{(\mu+\varepsilon)(\mu+\gamma)}\frac{\beta}{R_0} = \frac{\varepsilon}{(\mu+\varepsilon)(\mu+\gamma)}\frac{\beta(\mu+\varepsilon)(\mu+\gamma)}{\varepsilon\beta} = 1$$

The positive sign indicates there is a direct relationship between $\beta$ and $R_0$. The 1 means that, a unit increase in $\beta$, will result a unit increase in $R_0$.

For $\varepsilon$,

$$S_\varepsilon = \frac{\partial R_0}{\partial \varepsilon}\frac{\varepsilon}{R_0} = \frac{\partial}{\partial \varepsilon}\frac{\varepsilon\beta}{(\mu+\varepsilon)(\mu+v)}\frac{\varepsilon}{R_0} = \frac{(\mu+\varepsilon)\beta - \varepsilon\beta}{(\mu+\gamma)(\mu+\varepsilon)^2}\frac{\varepsilon(\mu+v)(\mu+\varepsilon)}{\varepsilon\beta}$$
$$= \frac{(\mu+\varepsilon)-\varepsilon}{(\mu+\varepsilon)} = \frac{\mu+\varepsilon}{\mu}$$

Which means an increase in $\varepsilon$ lead to a decrease in the basic reproduction number.

For $\mu$,

$$S_\mu = \frac{\partial R_0}{\partial \mu}\frac{\mu}{R_0} = \frac{\partial}{\partial \mu}\left(\frac{\varepsilon\beta}{(\mu+\varepsilon)(\mu+\gamma)}\right)\frac{\mu}{R_0} = -\frac{\varepsilon\beta(2\mu+\varepsilon+\gamma)}{(\mu+\varepsilon)^2(\mu+\gamma)^2}\frac{\mu(\mu+\varepsilon)(\mu+\gamma)}{\varepsilon\beta}$$
$$= -\frac{\mu(2\mu+\varepsilon+\gamma)}{(\mu+\gamma)(\mu+\varepsilon)}$$

This means that when $\mu$ get a decrease, then the basic reproduction number also get decreased.

For $\gamma$,

$$S_\gamma = \frac{\partial R_0}{\partial \gamma}\frac{\gamma}{R_0} = \frac{\partial}{\partial \gamma}\frac{\varepsilon\beta}{(\mu+\varepsilon)(\mu+v)}\frac{\gamma}{R_0} = \frac{-\varepsilon\beta}{(\mu+\varepsilon)}(\mu+\gamma)^{-2}\frac{\gamma}{R_0}$$
$$= \frac{-\varepsilon\beta\gamma(\mu+\varepsilon)(\mu+\gamma)}{(\mu+\varepsilon)(\mu+\gamma)^2}\cdot\varepsilon\beta = \frac{\gamma}{(\mu+\gamma)}$$

This implies that decreasing $\gamma$ lead to decrease the basic reproduction number $R_0$. The bar graph in Figure 5 shows the sensitivity index of the parameters influencing the basic reproduction number $R_0$.

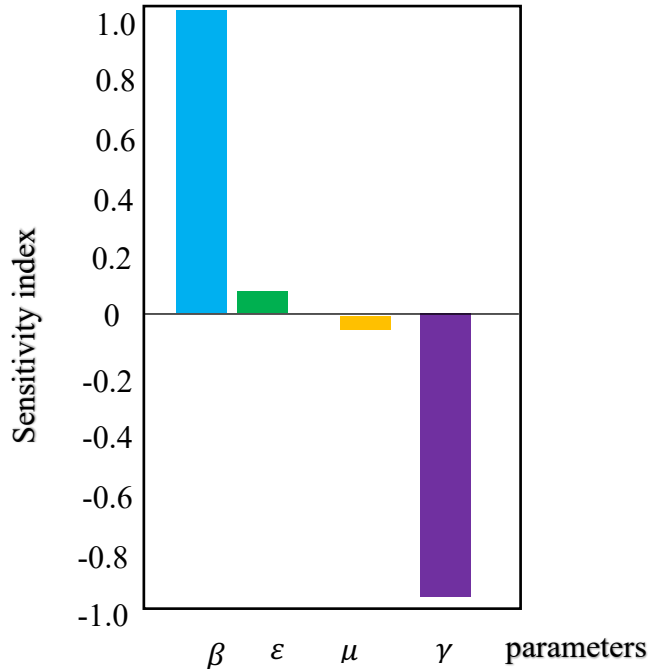

Figure 5: Sensitivity Index Bar

## 4.2: Numerical Simulations

In this section we are about to show numerical simulation of an epidemic model [31-36]. To simulate a system, we first consider a system.

$$\frac{dS}{dt} = -\beta SI/N$$
$$\frac{dE}{dt} = \beta SI/N - \varepsilon E$$
$$\frac{dI}{dt} = \varepsilon E - \gamma I$$
$$\frac{dR}{dt} = \gamma I$$

(4.1)

which is a classical SEIR model [8].
Where,

| Parameter | Identity |
|---|---|
| $\beta$ | Infection Rate |
| $\varepsilon$ | Rate at which infected become infectious |
| $\gamma$ | Recovery rate |
| N | Total population |
| S | Susceptible individuals |
| E | Infected but not infectious individuals |
| I | Infected and infectious individuals |
| R | Recovered individuals |

Table 4

To solve this, we assumed the parameter value as following:

| Parameter | Value | Data source |
|---|---|---|
| $\beta$ | 0.95 | assumed |
| $\varepsilon$ | 0.50 | assumed |
| $\gamma$ | 0.09 | assumed |
| N | 1000 | assumed |

Table 5

To solve the system, we also needed initial conditions and we also assumed the initial conditions which are shown in the following figure:

| Parameter | Value | Data source |
|---|---|---|
| $S_0$ | 960 | assumed |
| $E_0$ | 10 | assumed |
| $I_0$ | 30 | assumed |
| $R_0$ | 0 | assumed |

Table 6

We solved the system (4.1) using MATLAB and to solve this we used the RK-4 method to get numerical simulations using the above data given in the Table: 5 & 6 and we got the following simulations:

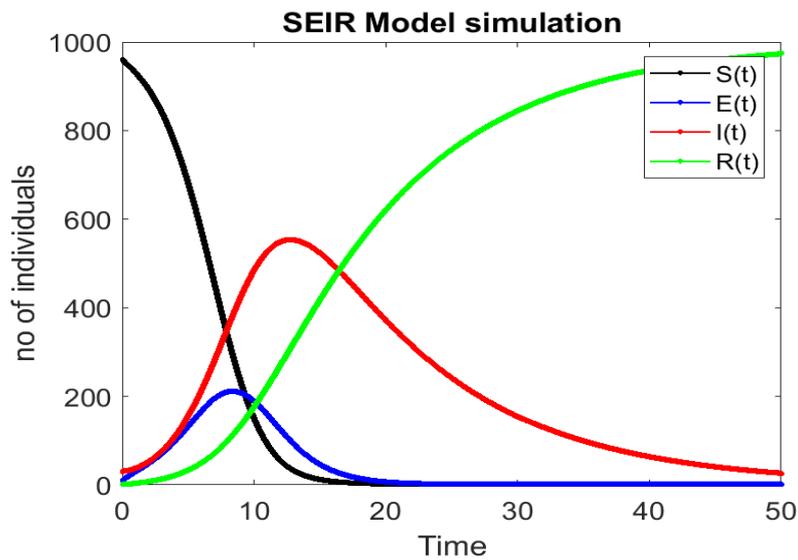

Figure 6: SEIR numerical simulation.

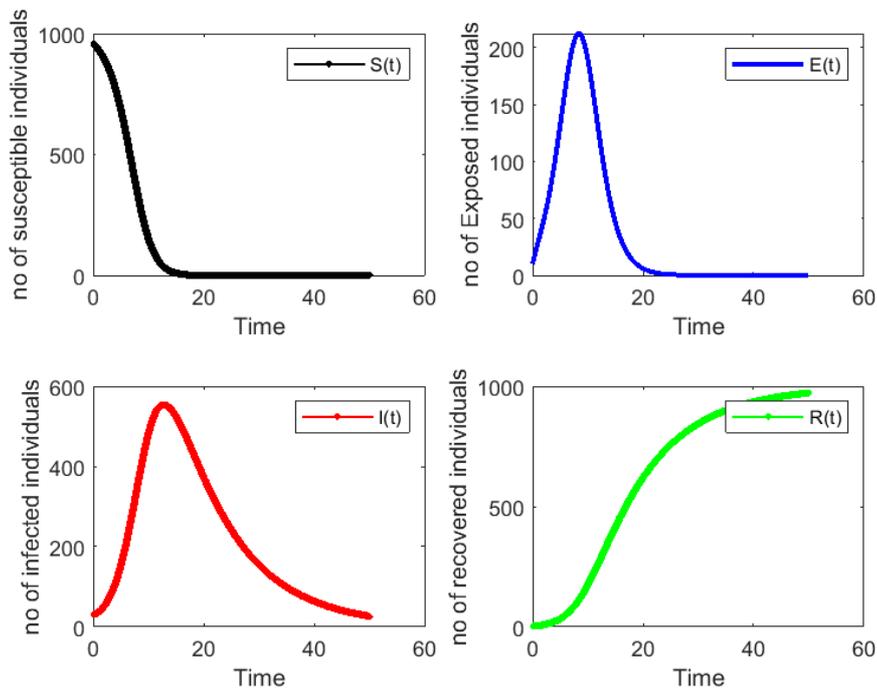

Figure 7: SEIR numerical simulation (separated)

In the above figure 6 & 7 the black line, blue line, red line, and the green line denoted the number of susceptible, infected but not infectious, infected with infectious and recovered individuals with respect to time [37-40]. We see from the above figure that initially the number of infected individuals increases with respect to time and after a certain period the infection rate goes down and as time tends to infinite the infection vanishes. And the number of

| Abbreviation | Full meaning |
|:---:|:---:|
| SEIR | Susceptible-exposed-infective-removed |
| DEE | Disease-free-equilibrium |
| EE | EE |
| ODES | Ordinary Differential Equations |
| RK4 | Fourth-order Runge–Kutta method |
| MATLAB | MATtrix LABoratory |

recovered individuals is increasing with respect to time. And it shows that as time tends to infinite the infection dies out.

## Conclusion:

The main objective of this work is to introduce all the concepts for studying mathematical modelling. On this project we have formulated a model for the population of any epidemiology. We divided the whole population into different classes: Susceptible, Exposed, Infected, and Recovered. The well-post-ness of the model has been established and the stability of the proposed model has been analyzed with the help of basic reproduction number [41-43]. In the system four types of equilibrium points have been discussed i.e. one is for disease free system and one is for endemic situation and the others are for one disease is present whenever the other is die out. Global stability has been analyzed by Lyapunov function. With the help of Bifurcation analysis, we have shown the qualitative changes in the behavior of different systems as parameters are varied. Numerical investigation has been presented by using MATLAB. At the end, by sensitive analysis we understood how changes in model parameters or initial conditions affect the outcomes of epidemic models.

However, diseases are produced by viruses that exist in various varieties, and a large percentage of human hosts that become infected are asymptotic, meaning they only show moderate symptoms, but still mobile. When administered appropriately, existing vaccines are crucial in preventing the disease's spread and averting the development of severe symptoms but diseases are always unpredictable. Therefore, this project aims to model any epidemiology and do all the theoretical analyses.

The dynamics of infectious disease can be explored precisely with the help of mathematical models. As most of our analyses are based on SEIR model, which perfectly describe the situation for the sake of liquidity. The diversity of SEIR model with effective conditions of the outbreak have been discussed with the help of some published model. Since the SEIR model deals with epidemic diseases which has latency period, a class of population that contains the exposed and asymptomatic individuals has been under consideration. These exposed population roam among susceptible population and transfers the pathogen unexpectedly. This work mainly focused on the nature of transmitting virus and changes of infection rate for the sake of exposed class. The impact of parameter from exposed to infected class population shows that the infection can be prevented by decreasing the value of this parameter, so the respective person can take appropriate action.